\documentclass[12pt]{article}
\usepackage{amsfonts,amsmath,amstext,amssymb}
\usepackage{dsfont}
\usepackage{graphicx}

\newtheorem{thm}{Theorem}
\newtheorem{lemma}[thm]{Lemma}
\newtheorem{corollary}{Corollary}
\newtheorem{remark}{Remark}

\newenvironment{demo}{\noindent \textit{Proof.} }{\hfill{$\Box$}}

\author{Alma L. Albujer\\
{\tt \small alma.albujer@uco.es}\\
\and 
Magdalena Caballero\\
{\tt \small magdalena.caballero@uco.es}\\
\and
{\small Departamento de Matem\'aticas, Campus de Rabanales,}\\ {\small Universidad de C\'ordoba, 14071 C\'ordoba, Spain}}

\date{}

\begin{document}

\title{Geometric properties of surfaces with the same mean curvature in $\mathds{R}^3$ and $\mathds{L}^3$}

\maketitle

\begin{abstract}
Spacelike surfaces in the Lorentz-Minkowski space $\mathds{L}^3$ can be endowed with two different Riemannian metrics, the metric inherited from $\mathds{L}^3$ and the one induced by the Euclidean metric of $\mathds{R}^3$. It is well known that the only surfaces with zero mean curvature with respect to both metrics are open pieces of the helicoid and of spacelike planes. We consider the general case of spacelike surfaces with the same mean curvature with respect to both metrics. One of our main results states that those surfaces have non-positive Gaussian curvature in $\mathds{R}^3$. As an application of this result, jointly with a general argument on the existence of elliptic points, we present several geometric consequences for the surfaces we are considering. Finally, as any spacelike surface in $\mathds{L}^3$ is locally a graph over a domain of the plane $x_3=0$, our surfaces are locally determined by the solutions to the $H_R=H_L$ surface equation. Some uniqueness results for the Dirichlet problem associated to this equation are given. 

\vspace{.3cm}


\end{abstract}

\section{Introduction}\label{intro}

A hypersurface in the Lorentz-Minkowski space $\mathds{L}^{n+1}$ is said to be spacelike if its induced metric is a Riemannian one. We can endow a spacelike hypersurface in $\mathds{L}^{n+1}$ with another Riemannian metric, the one inherited from the Euclidean space $\mathds{R}^{n+1}$. Therefore, we can consider two different mean curvature functions  on a spacelike hypersurface, the mean curvature function related to the metric induced by $\mathds{R}^{n+1}$, that we will denote by $H_R$, and the one related to the metric inherited from $\mathds{L}^{n+1}$, $H_L$. %In general both functions do not coincide.
 
A hypersurface in $\mathds{R}^{n+1}$ is said to be minimal if its mean curvature function vanishes identically, that is $H_R\equiv 0$. Analogously, a spacelike hypersurface in $\mathds{L}^{n+1}$ is said to be maximal if $H_L\equiv 0$. This terminology comes from the fact that minimal (maximal) hypersurfaces locally minimize (maximize) area among all nearby hypersurfaces sharing the same boundary, see \cite{Lopez}. 

The study of minimal and maximal hypersurfaces is a topic of wide interest. One of the main results about the global geometry of minimal surfaces is the well-known Bernstein theorem, proved by Bernstein~\cite{Bernstein} in 1915, which states that the only entire minimal graphs in $\mathds{R}^3$ are the planes. Some decades later, in 1970, Calabi~\cite{Calabi} proved its analogous version for spacelike surfaces in the Lorentz-Minkowski space, the Calabi-Bernstein theorem, which states that the only entire maximal graphs in $\mathds{L}^3$ are the spacelike planes. An important difference between both results is that the Bernstein theorem can be extended to minimal graphs in $\mathds{R}^{n+1}$ up to dimension $n=7$, but it is no longer true for bigger dimensions~\cite{BombieriDiGiorgiGiusti}. However, the Calabi-Bernstein theorem holds true for any dimension as it was proved by Calabi~\cite{Calabi} for dimension $n\leq 4$, and by Cheng and Yau~\cite{ChengYau} for arbitrary dimension. 

It is interesting to note that any complete spacelike hypersurface in $\mathds{L}^{n+1}$ is necessarily an entire graph over any spacelike hyperplane, see \cite[Proposition 3.3]{AliasRomeroSanchez}. Therefore, the Calabi-Bernstein theorem can also be expressed in a parametric way by asserting that the only complete maximal hypersurfaces in $\mathds{L}^{n+1}$ are the spacelike hyperplanes. This parametric version is not true in $\mathds{R}^{n+1}$, indeed there exists a wide family of examples of non-trivial complete minimal hypersurfaces in $\mathds{R}^{n+1}$. 

As an immediate consequence of the above results, we conclude that the only complete hypersurfaces that are simultaneously minimal in $\mathds{R}^{n+1}$ and maximal in $\mathds{L}^{n+1}$ are the spacelike hyperplanes. 

Going a step further, we can consider spacelike hypersurfaces with the same constant mean curvature functions $H_R$ and $H_L$. Heinz~\cite{Heinz}, Chern~\cite{Chern} and Flanders~\cite{Flanders} proved that the only entire graphs with constant mean curvature $H_R$ in $\mathds{R}^{n+1}$ are the minimal graphs. There are examples of entire spacelike graphs with constant mean curvature $H_L$ in $\mathds{L}^{n+1}$ which are not maximal, for instance the hyperbolic spaces. However, taking into account the Calabi-Bernstein theorem, we conclude again that the only complete spacelike hypersurfaces in $\mathds{L}^{n+1}$ with the same constant mean curvature functions $H_R$ and $H_L$ are the spacelike hyperplanes.

Without assuming any completeness hypothesis, Kobayashi~\cite{Kobayashi} studied the problem for $H_R=H_L=0$ in the $2$-dimensional case. He showed that the only surfaces that are simultaneously minimal and maximal are open pieces of a spacelike plane or of a helicoid in the region where the helicoid is spacelike. However, nothing is known neither for bigger dimension nor for more general mean curvature functions $H_R$ and $H_L$. 

Our purpose is to study some local and global geometric properties of the spacelike surfaces in $\mathds{L}^3$ such that $H_R=H_L$, not necessarily constant. Although we will focus on dimension 2, some results are still true in arbitrary dimension. 

It is well known that any spacelike surface can be locally seen as a graph over an open domain of a spacelike plane, which without loss of generality can be supposed to be the plane $x_3=0$, see \cite{Lopez}. That is, a spacelike surface is locally defined by a smooth function $u$. Therefore, the functions $H_R$ and $H_L$ can be written in terms of such a function $u$ and its partial derivatives. In this way, the identity $H_R=H_L$ becomes a quasilinear elliptic partial differential equation, everywhere except at those points at which the Euclidean gradient of $u$ vanishes, where the equation is parabolic. 

In Section~\ref{prelim} we present some basic preliminaries on spacelike hypersurfaces in $\mathds{L}^{n+1}$ and their mean curvature functions with respect to the metrics inherited from $\mathds{R}^{n+1}$ and $\mathds{L}^{n+1}$. In the next section we state a result on the existence of an elliptic point in a hypersurface of $\mathds{R}^{n+1}$ under some appropriate assumptions, and we see that the same result holds for spacelike hypersurfaces in $\mathds{L}^{n+1}$. 
In Section~\ref{nonconvex} we consider spacelike surfaces in $\mathds{L}^3$ such that $H_R=H_L$. We prove that for those surfaces $K_R$ is always non-positive, and if the mean curvature does not vanish at a point, then the surface is locally non-convex at that point, {\bf Theorem~\ref{thm:nonconvex}}. 
From this theorem, as well as from the result on the existence of an elliptic point, we get some consequences to which the rest of the section is devoted. Specifically, we prove the following theorems. 

\noindent{\bf Theorem~\ref{cor:convexhull}. }{\it Let $\Sigma$ be a compact spacelike surface with (necessarily) non-empty boundary such that $H_R=H_L$. Then $\Sigma$ is contained in the convex hull of its boundary.}

\noindent{\bf Theorem ~\ref{thm:asymptotic}. }{\it The only spacelike graphs $\Sigma_u$ in $\mathds{L}^3$ defined over a domain $\Omega\subseteq\mathds{R}^2$ of infinite width, with $H_R=H_L$, and asymptotic to a spacelike plane, are (pieces of) spacelike planes. }

\noindent{\bf Theorem~\ref{thm:diameter}. }{\it Let $\Sigma_u$ be a spacelike graph in $\mathds{L}^3$ over a domain $\Omega\subseteq\mathds{R}^2$ such that $H_R=H_L$. Then
$$
\mathrm{width}(\Omega^\ast)\leq\dfrac{1}{\sqrt{2}\,\inf_{\Omega^\ast}|H_L|},
$$ where $\Omega^\ast$ is the set of points at which the tangent plane is non-horizontal.}

Finally, in Section~\ref{consequences} we present the $H_R=H_L$ surface equation. Any spacelike surface is locally determined by a solution of this equation satisfying $|Du|<1$, where $D$ and $|\cdot|$ stand for the Euclidean gradient and Euclidean norm, respectively. We prove the uniqueness of the Dirichlet problem associated to this partial differential equation under some appropriate boundary conditions, {\bf Theorem~\ref{cor:dirichlet2}}. This is not trivial, since the equation is not always elliptic. 

\section{Preliminaries}\label{prelim}
Let $\mathds{L}^{n+1}$ be the $(n+1)$-dimensional Lorentz-Minkowski space, that is, $\mathds{R}^{n+1}$ endowed with the metric
\[
\langle,\rangle_L=dx_1^2+...+dx_{n}^2-dx_{n+1}^2,
\]
where $(x_1,...,x_{n+1})$ are the canonical coordinates in $\mathds{R}^{n+1}$, and let $|\cdot|_L$ denote its norm. It is easy to see that the Levi-Civita connections of the Euclidean space $\mathds{R}^{n+1}$ and the Lorentz-Minkowski space $\mathds{L}^{n+1}$ coincide, so we will just denote it by $\overline{\nabla}$.

A (connected) hypersurface $\Sigma^{n}$ in $\mathds{L}^{n+1}$ is said to be a spacelike hypersurface if $\mathds{L}^{n+1}$ induces a Riemannian metric on $\Sigma$, which is also denoted by $\langle,\rangle_L$. Given a spacelike hypersurface $\Sigma$, we can choose a unique future-directed unit normal vector field $N_L$ on $\Sigma$. %We will refer to $N_L$ as the future-pointing Gauss map of $\Sigma$. 
Let $\nabla^L$ denote the Levi-Civita connection in $\Sigma$ with respect to $\langle,\rangle_L$. Then the Gauss and Weingarten formulae for the spacelike hypersurface $\Sigma$ become
\begin{equation*} \label{eq:fGauss-l}
\overline{\nabla}_XY=\nabla^L_XY-\langle A_LX,Y\rangle_L N_L
\end{equation*}
and
\begin{equation*} \label{eq:fWein-l}
A_LX=-\overline{\nabla}_XN_L,
\end{equation*}
respectively, for any tangent vector fields $X,Y\in\mathfrak{X}(\Sigma)$, where $A_L:\mathfrak{X}(\Sigma)\rightarrow\mathfrak{X}(\Sigma)$ stands for the shape operator of $\Sigma$ with respect to $N_L$. The mean curvature function of $\Sigma$ with respect to $N_L$ is defined by
\[
H_L=-\frac{1}{n}\mathrm{tr}\,A_L=-\frac{1}{n}(k^L_1+...+k^L_{n}),\]
where $k^L_i,\,i=1,...,n$, stand for the principal curvatures of $(\Sigma,\langle,\rangle_L)$.

It is well known that there exists no closed (compact and without boundary) spacelike hypersurface in $\mathds{L}^{n+1}$~\cite{AledoAlias,AliasRomeroSanchez}. Therefore, every compact spacelike hypersurface $\Sigma$ in the Lorentz-Minkowski space necessarily has non-empty boundary.

The same topological hypersurface can also be considered as a hypersurface of the Euclidean space, that is $\mathds{R}^{n+1}$ with its usual Euclidean metric. For simplicity, we will just denote the Euclidean space by $\mathds{R}^{n+1}$, the Euclidean metric and the induced metric on $\Sigma$ by $\langle,\rangle_R$, and its norm by $|\cdot|_R$. In such a case, $\Sigma$ admits a unique upwards directed unit normal vector field, $N_R$. In an analogous way as in the Lorentzian case, let $\nabla^R$ denote the Levi-Civita connection in $\Sigma$ with respect to $\langle,\rangle_R$. The Gauss and Weingarten formulae read now 
\begin{equation} \label{eq:fGauss-r}
\overline{\nabla}_XY=\nabla^R_XY+\langle A_RX,Y\rangle_R N_R
\end{equation}
and
\begin{equation} \label{eq:fWein-r}
A_RX=-\overline{\nabla}_XN_R,
\end{equation}
respectively, $A_R:\mathfrak{X}(\Sigma)\rightarrow\mathfrak{X}(\Sigma)$ being the shape operator of $\Sigma$ with respect to $N_R$. The mean curvature function of $\Sigma$ with respect to $N_R$ is defined by
\[
H_R=\frac{1}{n}\mathrm{tr}\,A_R=\frac{1}{n}(k^R_1+...+k^R_{n}),\]
where $k^R_i,\,i=1,...,n$, stand for the principal curvatures of $(\Sigma,\langle,\rangle_R)$.

A spacelike hypersurface is locally a graph over an open domain of the hyperplane $x_{n+1}=0$, which can be identified with $\mathds{R}^n$. Therefore, for each $p\in\Sigma$ there exists an open neighborhood of $p$, $\Omega\subseteq\mathds{R}^n$, and a smooth function $u\in\mathcal{C}^\infty(\Omega)$ such that $\Sigma=\Sigma_u$ on this neighborhood, where
\[
\Sigma_u=\{(x_1,...,x_n,u(x_1,...,x_n)):(x_1,...,x_n)\in\Omega\}.
\]
It is easy to check that $\Sigma_u$ is a spacelike hypersurface if and only if $|Du|<1$, where $D$ and $|\cdot|$ stand for the gradient operator and the norm in the Euclidean space $\mathds{R}^n$, respectively. In this case, it is possible to get expressions for the normal vector fields $N_L$ and $N_R$, as well as for the mean curvature functions $H_L$ and $H_R$, in terms of $u$. Specifically, with a straightforward computation we get
\begin{equation}\label{eq:normLR}N_L=\dfrac{(Du,1)}{\sqrt{1- |Du|^2}} \qquad \textrm{and} \qquad N_R=\dfrac{(-Du,1)}{\sqrt{1+|Du|^2}}.\end{equation} And for the mean curvature functions we have
\begin{equation}\label{eq:HLR}
H_L = \dfrac{1}{2} \, \mathrm{div}\left(\dfrac{Du}{\sqrt{1- |Du|^2}}\right) \qquad \textrm{and} \qquad  
H_R = \dfrac{1}{2} \, \mathrm{div}\left(\dfrac{Du}{\sqrt{1+|Du|^2}}\right),
\end{equation}
where $\textrm{div}$ denotes the divergence operator in $\mathds{R}^n$. Let us observe that 
\[
\cos\theta=\frac{1}{\sqrt{1+|Du|^2}} \qquad \mathrm{and} \qquad \cosh\psi=\frac{1}{\sqrt{1-|Du|^2}},
\]
where $\theta$ and $\psi$ denote the angle between $N_R$ and $e_{n+1}=(0,...0,1)$ and the hyperbolic angle between $N_L$ and $e_{n+1}$, respectively. Moreover, from (\ref{eq:normLR}) it is immediate to get 
\begin{equation}\label{igualdad-con-angulos}
\dfrac{\langle X, N_L\rangle_L}{\cosh\psi}=-\dfrac{\langle X, N_R\rangle_R}{\cos\theta},
\end{equation}
for any $X\in\mathfrak{X}(\mathds{R}^{n+1})$ along $\Sigma$, which is a global equality since it does not depend on $u$. Let us observe that, in the previous expressions, we are writing $Du$ instead of $Du\circ\pi$, where $\pi$ is the canonical projection of $\Sigma_u$ onto $\Omega$. On behalf of simplicity, we will continue using this identification along the manuscript.  

\section{Setup}

As a first result, let us see that every compact hypersurface $\Sigma$ in $\mathds{R}^{n+1}$ with non-empty boundary, $\partial\Sigma$,  and such that it is not contained in the convex hull of $\partial\Sigma$, necessarily has an elliptic point, that is, a point at which all the principal curvatures have the same sign. The proof of this result follows the ideas of the proof of the classical result on the existence of an elliptic point in a compact surface without boundary in $\mathds{R}^3$. This idea was also used by the authors jointly with L\'{o}pez in~\cite{AlbujerCaballeroLopez} to prove the existence of an elliptic point in a non-planar compact spacelike surface in $\mathds{L}^3$ with planar boundary.

\begin{thm}
\label{lema:elliptic} Let $\Sigma$ be a compact hypersurface in $\mathds{R}^{n+1}$ with non-empty boundary. If $\Sigma$ is not contained in the convex hull of its boundary, $\textrm{conv}\,(\partial\Sigma)$, it necessarily has an elliptic point.
\end{thm} 

\begin{demo}
Let us assume that there exists a point $p'\in\Sigma\setminus\textrm{conv}\,(\partial\Sigma)$. Since $\textrm{conv}\,(\partial\Sigma)$ is a compact and convex set, by a classical result in convex geometry, $p'$ and $\textrm{conv}\,(\partial\Sigma)$ can be strongly separated. That is, there exists a hyperplane $\Pi$ in $\mathds{R}^{n+1}$ such that $p'$ and $\textrm{conv}\,(\partial\Sigma)$ are strictly contained in each one of the two open half-spaces determined by $\Pi$ (see for instance~\cite[Theorem 1.3.4]{Schneider}). Then, it is immediate to see that we can choose $q\in\mathds{R}^{n+1}$ and $\varrho>0$ such that $\textrm{conv}\,(\partial\Sigma)\subset B_q(\varrho)$ and $p'\in S_q(\varrho)$, where $B_q(\varrho)$ is the open ball in $\mathds{R}^{n+1}$ centered at $q$ of radius $\varrho$ and $S_q(\varrho)=\partial B_q(\varrho)$. Consequently, the function $f:\Sigma\rightarrow\mathds{R}$ defined by $f(p)=\langle p-q,p-q\rangle_R$ attains its maximum at a point $p^\ast\in \Sigma\setminus \textrm{conv}\,(\partial\Sigma)$, which is in particular an interior point, see Figure~\ref{fig:eliptico}.

\begin{figure}[h]
 \centerline{\includegraphics[width=6cm]{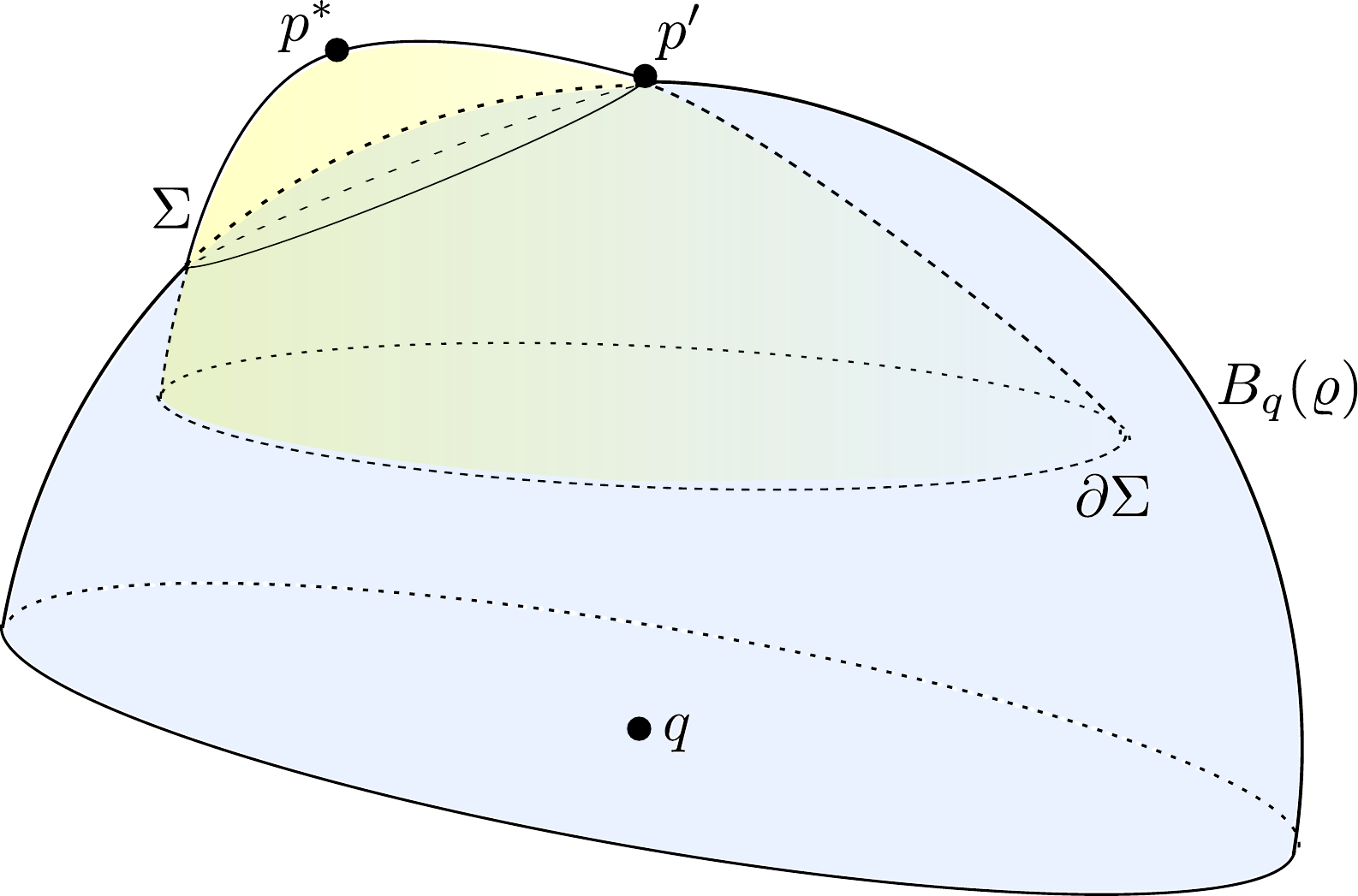}}  \caption{Existence of an elliptic point.}
 \label{fig:eliptico}
\end{figure}

Therefore, 
\begin{equation}
\label{eq:gradf} \nabla^R f(p^\ast)=0
\end{equation}
and
\begin{equation}
\label{eq:hessf} \mathrm{Hess}^Rf_{p^\ast}(v,v)\leq 0 \quad \forall v\in T_{p^\ast}\,\Sigma,
\end{equation}
where $\nabla^R$ and $\mathrm{Hess}^R$ denote the gradient and Hessian operators of $(\Sigma,\langle,\rangle_R)$.

Observe that for every $X\in\mathfrak{X}(\Sigma)$ and for every $p\in\Sigma$ it stands that
\[
X(p)(f)=2 \langle p-q, X(p) \rangle_R=\langle 2(p-q)^\top,X(p)\rangle_R,
\]
where $^\top$ denotes the tangent projection to $\Sigma$, which implies that
\[
\nabla^R f (p)= 2(p-q)^\top =2 (p-q-\langle p-q, N_R(p)\rangle_R N_R(p)).
\]
In particular, by \eqref{eq:gradf} we get
\begin{equation}
\label{eq:auxpast} \langle p^\ast-q,N_R(p^\ast)\rangle_R=\pm\varrho^\ast,
\end{equation}
where ${\varrho^\ast}^2=f(p^\ast)$.

On the other hand, by the Gauss \eqref{eq:fGauss-r} and Weingarten \eqref{eq:fWein-r} formulae we obtain
\[
\nabla^R_X\nabla^R f (p)=2(X(p)+\langle p-q, N_R(p) \rangle_R \, A_RX(p)).
\]
Therefore,
\begin{equation}
\label{eq:hessf2} \mathrm{Hess}^Rf_{p}(v,v)=2(\langle v,v\rangle_R+\langle p-q,N_R(p)\rangle_R \langle {A_R}_p(v),v\rangle_R), \quad \forall v\in T_p\Sigma.
\end{equation}

Let $\{e_1,...,e_n\}$ be an orthonormal basis of $T_{p^\ast}\Sigma$ which diagonalizes ${A_R}_{p^\ast}$, that is ${A_R}_{p^\ast}(e_i)=k^R_i(p^\ast)e_i$. Then \eqref{eq:hessf} and \eqref{eq:hessf2} yield
\[
\mathrm{Hess}^Rf_{p^\ast}(e_i,e_i)=2\left(1+\langle p^\ast-q,N_R(p^\ast)\rangle_R \, k^R_i(p^\ast)\right)\leq 0,\quad i=1,...,n,
\]
which jointly with \eqref{eq:auxpast} implies
\[
k^R_i(p^\ast)\leq -\frac{1}{\varrho^\ast}<0,\quad i=1,...,n, \quad \textrm{if} \quad \langle p^\ast-q,N_R(p^\ast)\rangle_R=\varrho^\ast >0,
\]
and
\[
k^R_i(p^\ast)\geq \frac{1}{\varrho^\ast}>0,\quad i=1,...,n, \quad \textrm{if} \quad \langle p^\ast-q,N_R(p^\ast)\rangle_R=-\varrho^\ast <0.
\]
Consequently, $p^\ast$ is an elliptic point.
\end{demo}

Our next goal is to show that Theorem~\ref{lema:elliptic} also holds for spacelike hypersurfaces in $\mathds{L}^{n+1}$. In fact, we will prove that if $\Sigma$ is  a spacelike hypersurface in $\mathds{L}^{n+1}$, a point in $\Sigma$ is elliptic with respect to the metric $\langle,\rangle_L$ if and only if it is elliptic with respect to $\langle,\rangle_R$, and so the result follows from Theorem~\ref{lema:elliptic}. To arrive to this conclusion we need the following previous result, which will also be used in the next section. 

\begin{lemma}\label{lema-curvaturas-normales}
Let $\Sigma$ be a spacelike hypersurface in $\mathds{L}^{n+1}$. Given $p\in\Sigma$ and $v\in T_p\Sigma$, let $\kappa_v^L(p)$ and $\kappa_v^R(p)$ denote the normal curvatures at $p$ in the direction of $v$ with respect to $\langle,\rangle_L$ and  $\langle,\rangle_R$, respectively. Then,
\[
\dfrac{|v|^2_R}{\cos\theta(p)}\,\kappa_v^R(p)=-\dfrac{|v|^2_L}{\cosh\psi(p)}\,\kappa_v^L(p).
\] 
\end{lemma}

\begin{demo}
Given $p\in\Sigma$ and $v\in T_p\Sigma$, let $\alpha$ be a smooth curve on $\Sigma$ such that $\alpha(0)=p$ and $\alpha'(0)=v$. We will work at $p$, but for simplicity we will omit it. Then, by definition,
\[
\kappa_v^R=k^R_\alpha\langle n_R,N_R\rangle_R \qquad \textrm{and} \qquad \kappa_v^L=-k^L_\alpha\langle n_L,N_L\rangle_L,
\]
where $k_\alpha^R$ and $k_\alpha^L$ denote the curvatures of $\alpha$ with respect to $\langle,\rangle_R$ and $\langle,\rangle_L$, respectively, while $n_R$ and $n_L$ are its normal vectors with respect to both metrics. Therefore, using the classical Frenet formulae we deduce
\begin{equation}\label{eq:curvnorm}
\kappa_v^R=\langle \overline{\nabla}_{t_R}t_R,N_R\rangle_R \qquad \textrm{and} \qquad \kappa_v^L=\langle \overline{\nabla}_{t_L}t_L,N_L\rangle_L,
\end{equation}
where $t_R=\frac{\alpha'}{|\alpha'|_R}$ and $t_L=\frac{\alpha'}{|\alpha'|_L}$. We combine~\eqref{igualdad-con-angulos} and~\eqref{eq:curvnorm} to finish the proof.
\end{demo}

Consequently, $\kappa_v^L$ and $\kappa_v^R$ always have opposite signs. And so, all the principal curvatures of $\Sigma$ with respect to $\langle,\rangle_R$ are positive if and only if all its principal curvatures with respect to $\langle,\rangle_L$ are negative, and vice versa. Hence, we have proved the Lorentzian version of Theorem~\ref{lema:elliptic}.

\begin{thm}
\label{cor:ellipticL} Let $\Sigma$ be a compact spacelike hypersurface in $\mathds{L}^{n+1}$ with necessarily non-empty boundary. If $\Sigma$ is not contained in the convex hull of its boundary, then it has an elliptic point.
\end{thm}

\section{Spacelike surfaces with $H_R=H_L$}\label{nonconvex}

In the following, we focus on dimension 2, and we deal with spacelike surfaces in $\mathds{L}^3$ endowed with the Riemannian metrics $\langle,\rangle_R$ and $\langle,\rangle_L$. The most important intrinsic curvature for surfaces is the well-known Gaussian curvature which is defined as \[K_R=\mathrm{det}(A_R)=k_1^Rk_2^R \qquad \textrm{and} \qquad K_L=-\mathrm{det}(A_L)=-k_1^Lk_2^L\] with respect to $\langle,\rangle_R$ and $\langle,\rangle_L$, respectively. Therefore, a point $p$ in a surface $\Sigma$ is elliptic with respect to the metric $\langle,\rangle_R$ if and only if $K_R(p)>0$, and with respect to $\langle,\rangle_L$ if and only if $K_L(p)<0$. In the case our surface is a graph over a domain of $\mathds{R}^2$, $\Sigma_u$, we can easily check that
\[
K_R=\frac{\mathrm{det}\,(D^2u)}{(1+|Du|^2)^2} \qquad \textrm{and} \qquad K_L=-\frac{\mathrm{det}\,(D^2u)}{(1-|Du|^2)^2},
\]
where $D^2u=\left(\begin{array}{cc}u_{xx}&u_{ xy}\\u_{xy}&u_{yy}\end{array}\right).$ From the above expressions it is clear that $K_RK_L\leq 0$ and $K_R(p)=0$, for any $p\in\Sigma$, if and only if $K_L(p)=0$. The last statement can also be deduced from Lemma~\ref{lema-curvaturas-normales}.

Going back to the mean curvature functions $H_R$ and $H_L$, it is interesting to observe that they have an expression in terms of the normal curvatures of any pair of orthogonal directions. Specifically,
\begin{equation}\label{eq:Hdirort}
H_R=\frac{1}{2}(\kappa^R_{v_1}+\kappa^R_{v_2}) \qquad \textrm{and} \qquad H_L=-\frac{1}{2}(\kappa^L_{w_1}+\kappa^L_{w_2}),
\end{equation}
where $\{v_1,v_2\}$ and $\{w_1,w_2\}$ are orthonormal basis of $T_p\Sigma$ with respect to $\langle,\rangle_R$ and $\langle,\rangle_L$, respectively.

\subsection{On the Gaussian curvature of surfaces with $H_R=H_L$}

\begin{thm}\label{thm:nonconvex}
Let $\Sigma$ be a spacelike surface in $\mathds{L}^3$ such that $H_R=H_L$. Then $K_R$ is everywhere non-positive. Moreover, if $K_R(p)=0$ for any $p\in\Sigma$, then $H_R(p)=0$. In particular, if the mean curvature at a point does not vanish, the surface is locally non-convex at that point. 
\end{thm}

\begin{demo}

We are going to work locally, so we can assume that there exist a domain $\Omega\subseteq\mathds{R}^2$ and a smooth function $u\in\mathcal{C}^\infty(\Omega)$ such that $\Sigma=\Sigma_u$. %Let us denote by $\pi$ its projection to $\mathds{R}^2$. 
We define $\Sigma^\ast$ as the graph of $u$ over the following open set
\begin{equation}\label{eq:Omegaast}
\Omega^\ast=\{(x_1,x_2)\in\Omega:Du(x_1,x_2)\neq 0\}.
\end{equation} 
Given $p\in\Sigma^\ast$, we consider its corresponding level curve contained in $\mathds{R}^2$ and we call its lifting to $\Sigma$, $\alpha$. We are working in a neighborhood of $p$, hence we can assume that $\alpha$ lies on $\Sigma^\ast$. Since $Du\neq 0$ in $\Omega^\ast$, its  distribution is integrable, so we can consider the integral curve through $\pi(p)$. We denote by  $\beta$ its lifting to $\Sigma^\ast$. Therefore, we have two curves defined on a neighborhood of $p$ which are orthogonal at $p$ for both $\langle,\rangle_R$ and $\langle,\rangle_L$.

The trace of the curve $\alpha$ is given by $\{(x,y,c):u(x,y)=c\}\subset\Sigma^\ast$ where $c=u(\pi(p))$, so a vector field tangent to $\alpha$ is $\alpha'=(-u_y,u_x,0)\circ\pi$. Then, Lemma~\ref{lema-curvaturas-normales} gives us the following relationship, where we have omitted the point $p$ on behalf of simplicity

\[
\kappa_{\alpha'}^R=-	\dfrac{|\alpha'|^2_L}{|\alpha'|^2_R}\dfrac{\cos\theta}{\cosh\psi}\,\kappa_{\alpha'}^L=-\sqrt{\frac{1-|Du|^2}{1+|Du|^2}}\,\kappa_{\alpha'}^L.
\]

On the other hand, it is easy to check that a tangent vector field to $\beta$ is given by $\beta'=(u_x,u_y,|Du|^2)\circ\pi$.  As a consequence of Lemma~\ref{lema-curvaturas-normales} we get

\[
\kappa_{\beta'}^R=-\dfrac{|\beta'|^2_L}{|\beta'|^2_R}\dfrac{\cos\theta}{\cosh\psi}\,\kappa_{\beta'}^L=-\left(\frac{1-|Du|^2}{1+|Du|^2}\right)^{3/2}\,\kappa_{\beta'}^L.
\]

By denoting $A=\sqrt{\frac{1-|Du|^2}{1+|Du|^2}}$, we rewrite the previous relations as 
\begin{equation}\label{eq:relkA1}
\kappa_{\alpha'}^R=-A\,\kappa_{\alpha'}^L \qquad \textrm{and} \qquad \kappa_{\beta'}^R=-A^3\,\kappa_{\beta'}^L.
\end{equation}

As we are dealing with two orthogonal curves at $p$ for both $\langle,\rangle_R$ and $\langle,\rangle_L$, and we are assuming  
$H_R=H_L$, from~\eqref{eq:Hdirort} we get
\begin{equation*}\label{eq:relkA2}
-\kappa_{\alpha'}^L-\kappa_{\beta'}^L=\kappa_{\alpha'}^R+\kappa_{\beta'}^R,
\end{equation*}
which jointly with \eqref{eq:relkA1} implies 
\begin{equation}\label{relacion-kLalfa-KLbeta}\kappa_{\alpha'}^L=-(A^2+A+1)\kappa_{\beta'}^L,
\end{equation}
and so
\[
\kappa_{\alpha'}^L \kappa_{\beta'}^L\leq 0.
\]
Consequently, the same inequality holds for the main curvatures of $\Sigma^\ast$, and so $K_L\geq 0$, or equivalently $K_R\leq 0$ in $\Sigma^\ast$. 

Consider now $p\in\Sigma\setminus\Sigma^\ast$. If $p\in\textrm{int}\,(\Sigma\setminus\Sigma^\ast)$, then $\Sigma$ is locally a horizontal plane, and so $K_L(p)=K_R(p)=0$. Otherwise $p\in\partial\Sigma^\ast$ and $K_R(p)\leq 0$ by a continuity argument. So the inequality $K_R\leq 0$ holds everywhere on $\Sigma$.

Finally, 
if we consider $p\in\Sigma^\ast$ such that $K_R(p)=0$, then $K_L(p)=0$ and there cannot exist normal curvatures at $p$ with opposite signs. So, from~\eqref{relacion-kLalfa-KLbeta} we obtain $\kappa_{\alpha'}^L(p)=\kappa_{\beta'}^L(p)=0$. Hence, $H_L(p)=0$. If $p\in\Sigma\setminus\Sigma^\ast$, we conclude the proof by using again a continuity argument.

\end{demo}

\subsection{Some consequences}

We have just proved that a surface $\Sigma$ such that $H_R=H_L$ does not have any elliptic point, which jointly which Theorem~\ref{lema:elliptic} lead to some interesting geometric consequences, to which the rest of the manuscript is devoted. The first of them is immediate from both results.

\begin{thm}\label{cor:convexhull}
Let $\Sigma$ be a compact spacelike surface with (necessarily) non-empty boundary such that $H_R=H_L$. Then $\Sigma$ is contained in the convex hull of its boundary.
\end{thm}

Let us recall that any spacelike surface is locally a graph $\Sigma_u$ over a domain $\Omega\subseteq\mathds{R}^2$. From now on, we will focus on spacelike graphs.

Firstly, we present a uniqueness result for graphs which are asymptotic to a spacelike plane, where the term asymptotic is defined as follows. We say that two entire graphs $\Sigma_u$ and $\Sigma_v$ over $\mathds{R}^2$ are \textit{asymptotic} if for every $\varepsilon>0$ there exists a compact set $K\subset\mathds{R}^2$ such that $|u(x,y)-v(x,y)|<\varepsilon$ for every $(x,y)\in \mathds{R}^2\setminus K$. Observe that, without loss of generality, we can consider that those compact sets are Euclidean discs of a certain radius. If we define the \textit{width} of a set in $\mathds{R}^2$ as the supremum of the diameter of the closed discs contained in it, the concept of asymptotic graphs is not only well-defined in the case of entire graphs, but also in the case of graphs over a domain of infinite width, that is, a domain which contains closed discs of any radius. Notice that this definition is a generalization of the classical concept of width for a convex body, see~\cite{Schneider}.

\begin{thm}\label{thm:asymptotic}
The only spacelike graphs $\Sigma_u$ in $\mathds{L}^3$ defined over a domain $\Omega\subseteq\mathds{R}^2$ of infinite width, with $H_R=H_L$, and asymptotic to a spacelike plane, are (pieces of) spacelike planes. 
\end{thm}

\begin{demo}
Let us notice that $\Sigma_u$ is a graph over any spacelike plane, and in particular, over the plane to which it is asymptotic. To prove it, it is enough to observe that if some timelike line intersects $\Sigma_u$ twice, the vertical plane containing that line cuts $\Sigma_u$ in a curve which is timelike at some point, which is a contradiction.

Let us denote by $\Pi$ the plane to which $\Sigma_u$ is asymptotic and let  $v\in \mathcal{C}^{\infty}(\Omega')$ be the function such that $\Sigma_u=\Sigma_v$, $\Omega'\subseteq\Pi$ being a domain. Notice that the width of $\Omega'$ is also infinite.

For any $\varepsilon >0$ there exists $(y_1,y_2)\in\Omega'$ and $R>0$ such that $|v(x_1,x_2)|<\varepsilon$ for every $(x_1,x_2)\in\Pi\setminus \bar{B}_{(y_1,y_2)}(R)$. By Theorem~\ref{cor:convexhull}, we know that the graph of the restriction of $v$ to $\bar{B}_{(y_1,y_2)}(R)$ is contained in the convex hull of its boundary. Therefore $|v(x_1,x_2)|\leq\varepsilon$ for all $(x_1,x_2)\in \bar{B}_{(y_1,y_2)}(R)$, so  this inequality holds globally on $\Omega'$. Taking limits when $\varepsilon$ tends to $0$, we conclude that $\Sigma_u=\Pi.$
\end{demo}

As another consequence, given a spacelike graph satisfying $H_R=H_L$, we can determine the width of $\Omega^\ast$, defined as in~\eqref{eq:Omegaast}, that is, the set of points on which the gradient of the function does not vanish, or equivalently, the points where the tangent plane to the surface is not a horizontal one. We need a previous result relating the mean curvature of those surfaces and the curvature of its level curves.

\begin{lemma}
Let $\Sigma_u$ be a spacelike graph in $\mathds{L}^3$ over a domain $\Omega\subseteq\mathds{R}^2$ such that $H_R=H_L$. If $\widetilde{k}$ denotes the curvature of the level curves of $u$, then for all $p\in\Sigma^\ast$  \begin{equation}\label{relacion-curvatura-R2-y-media}
|H_L(p)|<\frac{1}{2\sqrt{2}}|\widetilde{k}(\pi(p))|.  
\end{equation}
\end{lemma}

\begin{demo}
Following the notation introduced in the proof of Theorem~\ref{thm:nonconvex}, let $\alpha$ denote the lifting to $\Sigma^\ast$ of a level curve of $\Sigma$ through a point of $\Sigma^\ast$. We call $\widetilde{\alpha}$ such a level curve. It is possible to relate the normal curvature of $\Sigma^\ast$ along $\alpha$ in the direction of $\alpha'$ with the curvature of the planar curve $\widetilde{\alpha}$ in $\mathds{R}^2$. 
\[
\kappa_{\alpha'}^L=\langle\overline{\nabla}_{t_L}t_L,N_L\rangle_L=-\frac{|Du|}{\sqrt{1-|Du|^2}}\langle D_{\widetilde{t}}\,\widetilde{t},\widetilde{n}\rangle_{\mathds{R}^2}\circ\pi=-\widetilde{k}\circ\pi\frac{|Du|}{\sqrt{1-|Du|^2}},
\]
where $D$ and $\langle,\rangle_{\mathds{R}^2}$ stand for the Levi-Civita connection and the usual metric of the Euclidean plane $\mathds{R}^2$, respectively, and $\widetilde{t}=\dfrac{(-u_y,u_x)}{|Du|},\widetilde{n}=\dfrac{-Du}{|Du|},\widetilde{k}$ is the Frenet apparatus of $\widetilde{\alpha}$ as a planar curve of $\mathds{R}^2$. Therefore, from~\eqref{relacion-kLalfa-KLbeta} we get
\[
2H_L=\frac{A+1}{A^2+A+1}\frac{|Du|}{\sqrt{1+|Du|^2}}\widetilde{k}\circ\pi=f(|Du|){\widetilde{k}}\circ\pi.
\]
Since $f(|Du|)$ is positive, $f$ is increasing in $|Du|$ and $|Du|<1$, we get~\eqref{relacion-curvatura-R2-y-media}.

\end{demo}

\begin{remark}
It is interesting to notice that the non-strict inequality holds for any point at which a level curve is defined. The reason is that those points belong to $\overline{\Omega^\ast}$.
\end{remark}

\begin{thm}\label{thm:diameter}
Let $\Sigma_u$ be a spacelike graph in $\mathds{L}^3$ over a domain $\Omega\subseteq\mathds{R}^2$ such that $H_R=H_L$. Then
\begin{equation}\label{width-inequality}
\mathrm{width}(\Omega^\ast)\leq\dfrac{1}{\sqrt{2}\,\inf_{\Omega^\ast}|H_L|}.
\end{equation} 
\end{thm}

\begin{demo}

If $\inf_{\Omega^\ast}|H_L|=0$, there is nothing to prove. 

Otherwise, we have $|H_L|\geq \inf_{\Omega^\ast}|H_L| =c>0$ in $\Sigma^\ast$.
And, as a consequence of \eqref{relacion-curvatura-R2-y-media}, we get
\begin{equation}\label{eq:descurv}
|\widetilde{\kappa}|>2\sqrt{2}\,c>0 \quad \textrm{in} \quad \Omega^\ast.
\end{equation} 

Firs of all, let us notice that $\Omega^\ast$ is an open set of $\mathds{R}^2$. We consider all the level curves in $\Omega^\ast$, we order them by the value of $u$ on each of them and we orient them in a way such that its normal vectors point to the direction on which $u$ decreases. 

We proceed by reductio ad absurdum assuming that the width of $\Omega^\ast$ is bigger than $\frac{1}{\sqrt{2}\,c}$. Then, there exists a point $q\in\Omega^\ast$ such that $\bar{B}_q\left(1/2\sqrt{2}\,c\right)\subset\textrm{int}\, \Omega^\ast$. Since $\bar{B}_q\left(1/2\sqrt{2}\,c\right)$ is compact, $u$ attains a maximum in it. Even more, $Du$ does not vanish in $B_q\left(1/2\sqrt{2}\,c\right)$, and so this extremal value is only attained on the boundary of the disc. 

We pick a point $p$ at which a maximum is attained. The level curve through $p$ lie in $\Omega^\ast\setminus B_q\left(1/2\sqrt{2}\,c\right)$. And so, it is tangent to the boundary of the disc at $p$. The normal vector to the curve at $p$ points to the interior of the disc, while the curve is not locally contained in it, see Figure~\ref{fig:disco}. Consequently, \eqref{eq:descurv} does not hold at $p$, which is a contradiction.

\begin{figure}[h]
 \centerline{\includegraphics[width=6cm]{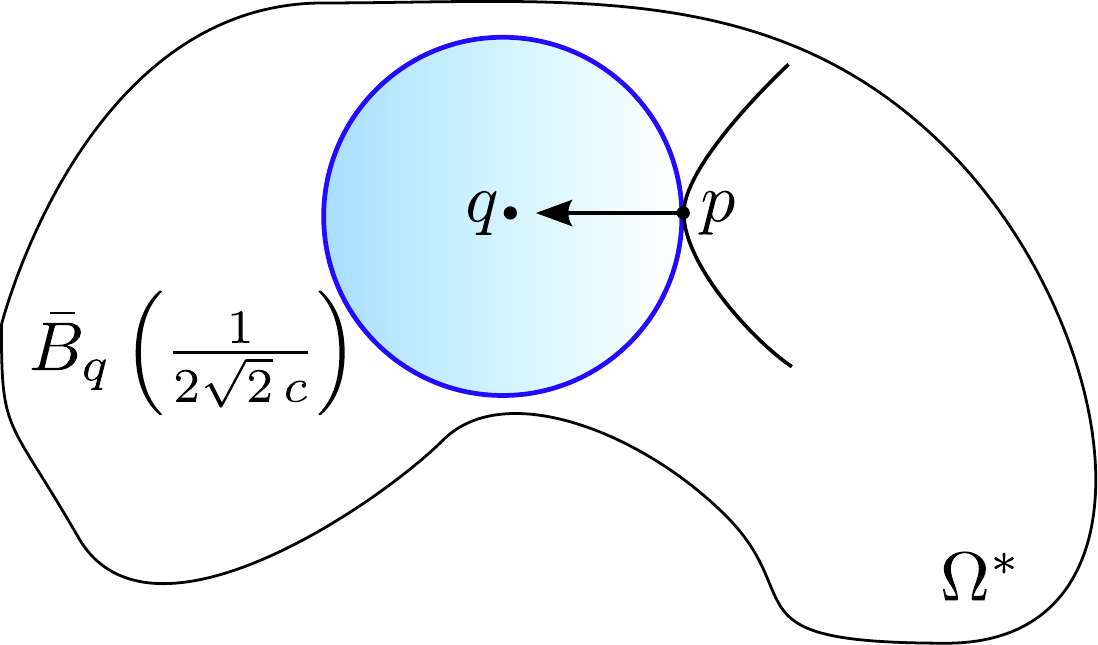}}  \caption{Level curve at $p$.} 
 \label{fig:disco}
\end{figure}

\end{demo}

\begin{remark}
In 1955 Heinz~\cite{Heinz} proved that if $\Sigma_u$ is a graph over a disc of radius $R$ and $|H_R|\geq c>0$, then $R\leq \frac{1}{c}$. When applied to a surface $\Sigma_u$ over a domain $\Omega$ such that $H_R=H_L$ and $|H_R|\geq c>0$, it tells us that $\mathrm{width}(\Omega^\ast)\leq\dfrac{2}{c}$. Inequality~\eqref{width-inequality} gives us an improved bound.
\end{remark}

As a direct consequence of Theorem~\ref{thm:diameter}, we get the following result.

\begin{corollary}\label{cor:diameter}
Let $\Sigma_u$ be a spacelike graph in $\mathds{L}^3$ such that $H_R=H_L$ and assume that $\Omega^\ast$ is a set of infinite width. Then $\inf_{\Sigma_u}|H_L|=0$.

Equivalently, there do not exist spacelike graphs satisfying $H_R=H_L$, $|H_L|\geq c$ for a certain constant $c>0$ and $\mathrm{width}(\Omega^\ast)=\infty$. 
\end{corollary}

\section{A quasi-linear PDE related to spacelike surfaces with $H_R=H_L$}\label{consequences}

Once again, we recall that any spacelike surface is locally a graph $\Sigma_u$ over a domain $\Omega\subseteq\mathds{R}^2$. From Theorem~\ref{thm:asymptotic} till the end of last section, we have focused on spacelike graphs satisfying $H_R=H_L$. Thanks to \eqref{eq:HLR}, if we consider the differential operator given by
\[
Q(u)=\mathrm{div}\left(\dfrac{Du}{\sqrt{1- |Du|^2}}\right) -   \mathrm{div}\left(\dfrac{Du}{\sqrt{1+|Du|^2}}\right),
\] 
those graphs are the solutions to the equation
\begin{equation}\label{eq:pde}
Q(u)=0,
\end{equation}
satisfying $|Du|<1$. We will refer to the above equation as {\it the $H_R=H_L$ surface equation}. 

Let us observe firstly that~\eqref{eq:pde} is an elliptic quasi-linear partial differential equation, everywhere except at those points where $Du=0$, at which it is parabolic. In fact, it is a straightforward computation to show that~\eqref{eq:pde} can be expressed as
\[
a_{11}u_{x_1x_1}+2a_{12}u_{x_1x_2}+a_{22}u_{x_2x_2}=0,
\]
$a_{11}$, $a_{12}$ and $a_{22}$ being smooth functions on $Du$ such that
\[
a_{11}a_{22}-a_{12}^2=\left(\frac{1}{\sqrt{1-|Du|^2}}-\frac{1}{\sqrt{1+|Du|^2}}\right)^2 \geq 0,
\]
with equality if and only if $Du=0$.

By a well-known result by Bernstein~\cite{Bernstein2}, the solutions of an elliptic quasi-linear partial differential equation for an analytic operator $Q$ are always analytic. Therefore, if $u$ is a solution of~\eqref{eq:pde} with $0<|Du|<1$, it is necessarily analytic. However, in general the analyticity of the solutions of~\eqref{eq:pde} cannot be guaranteed.

Let $\Omega\subset\mathds{R}^2$ be a bounded domain and $\psi\in\mathcal{C}^0(\partial\Omega)$. The Dirichlet problem related to the $H_R=H_L$ surface equation consists in finding a solution $u\in\mathcal{C}^2(\Omega)\cap\mathcal{C}^0(\overline{\Omega})$ to the boundary value problem
\begin{equation}\label{eq:Dirichlet}
\left.\begin{array}{ccc}
Q(u)=0 & \mathrm{in} & \Omega\\
|Du|<1 & \mathrm{in} & \Omega\\
u=\psi & \mathrm{on} & \partial\Omega
\end{array}\right\}.
\end{equation}

As a consequence of a uniqueness theorem for the Dirichlet problem associated to quasilinear elliptic operators ~\cite[Theorem 9.3]{GilbargTrudinger}, we get our next result.

\begin{thm}
Let $\Omega\subset\mathds{R}^2$ be a bounded domain with smooth boundary and $\psi\in\mathcal{C}^0(\partial\Omega)$ such that the Dirichlet problem~\eqref{eq:Dirichlet} admits a solution $u$ without critical points. Then, the solution is unique.
\end{thm}

\begin{remark} It is interesting to observe that~\cite[Theorem 9.3]{GilbargTrudinger} holds under four assumptions on the operator defining the equation, one of which does not hold in our case. Since we are assuming the spatially condition $|Du|<1$, the coefficients of $Q$ are not well-defined on the whole $\Omega\times\mathds{R}\times\mathds{R}^2$, as it is required in the last hypothesis, but just on $\Omega\times\mathds{R}\times B_0(1)$. However, studying in detail the proof of the cited theorem, we can realize that it is sufficient to consider the coefficients defined on $\Omega\times\mathds{R}\times B_0(1)$. 
\end{remark}

It is still more interesting to put emphasis on the fact that the proof does not work if the ellipticity fails somewhere. Therefore, we can not omit the hypothesis on the gradient of $u$. However, as a consequence of Theorem~\ref{cor:convexhull}, we get the following result on the uniqueness of the Dirichlet problem under appropriate boundary values.

\begin{thm}\label{cor:dirichlet2}
The only solutions to the Dirichlet problem~\eqref{eq:Dirichlet} with affine boundary values are the affine functions.
\end{thm}
\begin{demo}
Let $u$ be a solution of~\eqref{eq:Dirichlet} and $\Sigma_u$ its associated graph. From Theorem~\ref{cor:convexhull}, $\Sigma_u$ is contained in the convex hull of the graph of $\psi$, which is a planar curve. Therefore, the spacelike graph $\Sigma_u$ must be contained in a plane, and consequently $u$ is affine.

In the previous reasoning it is crucial to observe that Theorem~\ref{cor:convexhull} also works for $\mathcal{C}^2$- surfaces.

\end{demo}

\section*{Acknowledgements}

The first author is partially supported by MINECO/FEDER project reference MTM2015-65430-P, Spain, and Fundaci\'on S\'eneca project
reference 19901/ GERM/15, Spain. Her work is a result of the activity developed within the framework of the Program in Support of Excellence Groups of the
Regi\'{o}n de Murcia, Spain, by Fundaci\'{o}n S\'{e}neca, Science and Technology Agency of the Regi\'{o}n de Murcia.
The second author is partially supported by the Spanish Ministry of Economy and Competitiveness and European Regional Development Fund (ERDF), project MTM2013-47828-C2-1-P.

\bibliographystyle{amsplain}

\end{document}